\documentclass[12pt,reqno]{amsart}
\usepackage{amsmath,amssymb,amsfonts,amscd,latexsym,amsthm,mathrsfs,verbatim}
\usepackage[unicode]{hyperref}
\usepackage{fouridx}

\newtheorem*{mtheorem}{Main theorem}

\renewcommand{\d}{{\mathrm d}}
\newcommand{\m}{{\mathrm m}}

\begin{document}

\title{Ap\'ery limits and Mahler measures}

\author{Wadim Zudilin}
\address{Department of Mathematics, IMAPP, Radboud University, PO Box 9010, 6500~GL Nij\-megen, Netherlands}
\email{w.zudilin@math.ru.nl}

\dedicatory{To the blessed memory of Asmus L.~Schmidt}

\date{27 September 2021}

\begin{abstract}
It is the first paper which relates Ap\'ery limits to Mahler measures.
\end{abstract}


\maketitle

\section{Introduction}
\label{s1}

There seems to be no paper yet in which the two notions in the title, \emph{Ap\'ery limits} and \emph{Mahler measures}, interplay.
Though we marry them here somewhat artificially, there is a story behind the marriage.

\begin{mtheorem}
Consider the second order polynomial recursion
\begin{align*}
&
32(n+1)(2n+1)(2n+3)(4n+3)(4n+5)(8n+3)^2(8n+5)^2(8n+7)^2(8n+9)^2
\\ &\;\quad\times
(654311424n^8-763363328n^7+336592896n^6-62390272n^5+1949696n^4
\\ &\;\qquad
+706560n^3-43520n^2-2880n+225) \, r_{n+1}
\displaybreak[2]\\ &\;
-3(2n+1)(1011688620292508221440n^{20}+6913205571998806179840n^{19}
\\ &\;\qquad
+20264522984105850175488n^{18}+32694809236419354034176n^{17}
\\ &\;\qquad
+30311863475685170348032n^{16}+13940981944461731823616n^{15}
\\ &\;\qquad
-511604730729009774592n^{14}-3973361701783984930816n^{13}
\\ &\;\qquad
-1526200168532215332864n^{12}+194737260750799110144n^{11}
\\ &\;\qquad
+269867025941865168896n^{10}+37629485093249613824n^9
\\ &\;\qquad
-17598923437928087552n^8-5481696915139592192n^7+199010711963172864n^6
\\ &\;\qquad
+237194722753118208n^5+18366301200549888n^4-2095780639795200n^3
\\ &\;\qquad
-198344332843200n^2+7359342480000n+776998726875) \, r_n
\displaybreak[2]\\ &\;
-1536n^3(6n-1)(6n+1)(8n+1)^2(8n-1)^2(8n-3)^2(8n-5)^2
\\ &\;\quad\times
(654311424n^8+4471128064n^7+13313769472n^6+22567976960n^5+23822974976n^4
\\ &\;\qquad
+16040183808n^3+6728855040n^2+1608382656n+167760801) \, r_{n-1}
=0
\end{align*}
and its two solutions $\{q_n\}_{n=0}^\infty$, $\{p_n\}_{n=0}^\infty$ defined through the \textup(rational\textup!\textup) initial values
$$
q_0=p_0=1, \qquad q_1=\frac{1289}{160}, \;\; p_1=\frac{136185509}{15876000}.
$$
Then
\begin{equation}
\lim_{n\to\infty}\frac{p_n}{q_n}
=-\frac1{16\pi\sqrt2}\iint\limits_{|x|=|y|=1}\log|(x^4+1)y^2-2(x^4-4x^2+1)y+x^4+1|\,\frac{\d x}{x}\,\frac{\d y}{y}.
\label{A=MM}
\end{equation}
\end{mtheorem}

To clarify the matters, limits of quotients of solutions to linear recurrence equations, like the one featured in the left-hand side of \eqref{A=MM}, are usually called Ap\'ery limits \cite{ASZ08,CS20} (see also Section~\ref{s5} below).%
\footnote{The notion is `canonical' for the second order recurrence equations (with coefficients from $\mathbb Q[n]$), because there are two linear independent rational-valued solutions in this case and the quotient is well defined up to a linear-fractional transformation.
Here we deal with such second order recursions.}
The double integral on the right-hand side in \eqref{A=MM} represents (up to a multiple) the (logarithmic) Mahler measure \cite{Bo98,BZ20}
\begin{align*}
\m(P(x,y))
&=\frac1{(2\pi i)^2}\iint\limits_{|x|=|y|=1}\log|P(x,y)|\,\frac{\d x}{x}\,\frac{\d y}{y}
\end{align*}
of a particular polynomial $P(x,y)\in\mathbb Z[x,y]$.

\section{Manifestation}
\label{s2}

In our proof of the main theorem we simply show that the two sides in \eqref{A=MM} represent the \emph{same} quantity.
This number is somewhat exceptional (as it is one of a few that can be realised simultaneously as the Ap\'ery limit of a second order recursion and as a Mahler measure),
so it serves as the main hero of the story.

\section{The left-hand side}
\label{s3}

During this author's visit in the University of Copenhagen in 2004, A.~Schmidt suggested to look for second order polynomial recursions giving reasonable approximations to the following Dirichlet $L$-series at the point $s=2$:
$$
L(\chi_{-8},s)=\sum_{n=1}^\infty\frac{\chi_{-8}(n)}{n^s},
$$
where the (odd) quadratic character $\chi_{-8}=\big(\frac{-8}{\cdot}\big)$ modulo 8 is determined by the data $\chi_{-8}(1)=\chi_{-8}(3)=1$, $\chi_{-8}(5)=\chi_{-8}(7)=-1$.
The value $L(\chi_{-8},2)$ appears quite natural in a problem of approximations of complex
numbers by elements from $\mathbb Z[\sqrt{-2}]$ (see \cite{Pe33,Sc11}; the $L$-value is not featured there but its appearance is clarified in \cite[Theorem 3.1]{PP14}).
In the corresponding problems of approximations of real numbers by rationals and of complex numbers by elements from $\mathbb Z[\sqrt{-1}]$, the (`ergodic') analogues of $L(\chi_{-8},2)$ are $\pi^2/12$ and $G=L(\chi_{-4},2)$ (Catalan's constant).%
\footnote{Schmidt famously developed generalisations of the theory of continued fractions to the complex numbers \cite{Sc75,Sc82} and gave numerous explicit results for diophantine approximations in imaginary quadratic fields.}
The two latter constants admit second order Ap\'ery-like recursions; there were several constructions at that time known, which appeared in print later \cite{Za09,Zu05}.
But $L(\chi_{-8},2)$ looked quite challenging, and it is surprising to see that it did not show up as such an Ap\'ery limit since then.

The Ap\'ery-limit realisation of $L(\chi_{-8},2)$ comes from mixing the hypergeometric constructions~\cite{Zu02,Zu05}.
Executing Zeilberger's creative telescoping on the series
\[
r_n=\sum_{\nu=1}^\infty\frac{\d}{\d t}
\bigg(\frac{(-1)^{n+1}(2t+n-\frac12)\prod_{j=1}^n(t-j)^2(t+n+j-\frac12)^2}{2^9\prod_{l=0}^n(t+l+\frac18)(t+l-\frac18)(t+l-\frac38)(t+l-\frac58)}\bigg)\bigg|_{t=\nu}
\]
we find out that $r_n=q_nL(\chi_{-8},2)-p_n$, together with the coefficient sequences $q_n$ and $p_n$, satisfy the the linear homogeneous difference equation displayed in the main theorem.
The characteristic polynomial of the equation is
\[
\lambda^2-270\lambda-27=\big(\lambda-(135+78\sqrt3)\big)\big(\lambda-(135-78\sqrt3)\big),
\]
and a standard analysis leads to the asymptotics
\begin{equation}
\lim_{n\to\infty}q_n^{1/n}=\lim_{n\to\infty}p_n^{1/n}=78\sqrt3+135,
\quad
\lim_{n\to\infty}|r_n|^{1/n}=78\sqrt3-135=0.09996299\dotsc.
\label{asym}
\end{equation}
In particular, this demonstrates that $p_n/q_n\to L(\chi_{-8},2)$ as $n\to\infty$ and decodes the left-hand side of~\eqref{A=MM}.

In general, given (rational) $0<\alpha,\beta,\gamma<1$ with $\alpha\ne\beta$ and $\alpha+\beta\ne\gamma$, the choice
\[
R_n(t)=\frac{(-1)^{n+1}(2t+n-\gamma)\prod_{j=1}^n(t-j)^2(t+n+j-\gamma)^2}
{\prod_{l=0}^n(t+l-\alpha)(t+l-\beta)(t+l+\alpha-\gamma)(t+l+\beta-\gamma)}
\]
and launch of Zeilberger's creative telescoping
result in a second order recursion (whose coefficients are polynomials of degree~23 in~$n$) for the quantities
\[
r_n=\sum_{\nu=1}^\infty\frac{\d R_n(t)}{\d t}\bigg|_{t=\nu}
=q_nC(\alpha,\beta,\gamma)-p_n
\]
approximating the Ap\'ery limit
$$
C(\alpha,\beta,\gamma)
=\psi_1(1-\alpha)-\psi_1(1-\beta)+\psi_1(1+\alpha-\gamma)-\psi_1(1+\beta-\gamma),
$$
where
$$
\psi_1(x)=\frac{\d^2}{\d x^2}\log\Gamma(x)
=\sum_{n=0}^\infty\frac1{(n+x)^2}
$$
is the so-called trigamma function.
Notice that the characteristic polynomial of the recursion and asymptotics \eqref{asym} remain independent of the data $\alpha,\beta,\gamma$.
Since
$$
\psi_1(x)+\psi_1(1-x)=\frac{\pi^2}{\sin^2\pi x},
$$
the quantity
$$
C(\alpha,\beta,1)
=\pi^2\bigg(\frac1{\sin^2\pi\alpha}-\frac1{\sin^2\pi\beta}\bigg)
$$
is an algebraic multiple of $\pi^2$; here are examples of such Ap\'ery limits when $\gamma=1$:
\begin{gather*}
\begin{alignedat}{2}
\frac14C\bigg(\frac18,\frac38,1\bigg)&=\pi^2\sqrt2,
&\qquad
\frac18C\bigg(\frac1{12},\frac5{12},1\bigg)&=\pi^2\sqrt3,
\\
\frac54C\bigg(\frac15,\frac25,1\bigg)&=\pi^2\sqrt5,
&\qquad
\frac14C\bigg(\frac1{10},\frac3{10},1\bigg)&=\pi^2\sqrt5,
\end{alignedat}
\\
\frac18C\bigg(\frac19,\frac49,1\bigg)=\pi^2\cos\frac{\pi}9,
\quad
\frac18C\bigg(\frac19,\frac29,1\bigg)=\pi^2\cos\frac{2\pi}9,
\quad
\frac18C\bigg(\frac29,\frac49,1\bigg)=\pi^2\cos\frac{4\pi}9,
\\
\begin{alignedat}{2}
\frac14C\bigg(\frac1{15},\frac4{15},1\bigg)&=\pi^2\textstyle\sqrt3\sqrt{5+2\sqrt5},
&\quad
\frac14C\bigg(\frac2{15},\frac7{15},1\bigg)&=\pi^2\textstyle\sqrt3\sqrt{5-2\sqrt5},
\\
\frac1{12}C\bigg(\frac1{10},\frac25,1\bigg)&=\pi^2\bigg(\frac13+\frac1{\sqrt5}\bigg),
&\quad
\frac18C\bigg(\frac1{10},\frac15,1\bigg)&=\pi^2\bigg(\frac12+\frac1{\sqrt5}\bigg),
\\
\frac1{12}C\bigg(\frac3{10},\frac15,1\bigg)&=\pi^2\bigg(\frac13-\frac1{\sqrt5}\bigg),
&\quad
\frac18C\bigg(\frac3{10},\frac25,1\bigg)&=\pi^2\bigg(\frac12-\frac1{\sqrt5}\bigg).
\end{alignedat}
\end{gather*}
When $\gamma\ne1$, we find only two instances that can be identified with Dirichlet $L$-values (for odd characters):
\[
\frac1{160}C\bigg(-\frac1{12},\frac1{12},\frac12\bigg)-\frac1{160}
=L(\chi_{-4},2),
\quad
\frac1{64}C\bigg(-\frac18,\frac18,\frac12\bigg)-\frac1{64}
=L(\chi_{-8},2).
\]
Numerous realisations of Catalan's constant $L(\chi_{-4},2)$ as the Ap\'ery limits of second order recursions are already known \cite{Za09,Zu02}, while the coverage of $L(\chi_{-8},2)$ explicitly given in the main theorem is new.
In a certain sense, this explains the uniqueness of the latter quantity.

\section{The right-hand side}
\label{s4}

The Mahler measure evaluation
\[
\m\big((x^4+1)y^2-2(x^4-4x^2+1)y+x^4+1\big)
=L'(\chi_{-8},-1)
\]
for the right-hand side in \eqref{A=MM} is borrowed from the paper by G.~Ray \cite[Proposition 2]{Ra87}. The connection with $L(\chi_{-8},2)$ comes from the general formula
\[
L(\chi_{-N},2)=\frac{4\pi}{N\sqrt N}\,L'(\chi_{-N},-1),
\]
which is a consequence of the functional equation for $L(\chi_{-N},s)$.
Ray's formula can be replaced with the evaluations
\begin{align*}
\m\big((x+1)^4y-(x-1)^2(x^2+1)\big) &= L'(\chi_{-8},-1),
\\
\m\big((x+1)^2(x^2+x+1)y-(x^2+1)^2\big) &= \frac23 L'(\chi_{-8},-1),
\\
\m\big((x+1)^{12}y-(x-1)^8(x^4-x^3+x^2-x+1)\big) &= \frac{16}5 L'(\chi_{-8},-1)
\\ \intertext{given by D.~Boyd and F.~Rodriguez-Villegas in \cite{BV02}, or}
\m\big((x^2+1)y+(x+1)^2(x^2-x+1)\big) &= \frac13 L'(\chi_{-8}, -1)
\end{align*}
of F.~Brunault~\cite{Br21}.
An interest in casting the $L$-values $L'(\chi_{-N},-1)$ corresponding to odd Dirichlet characters $\chi_{-N}$ as two-variable Mahler measures has its roots in the first such example
\[
\m(x+y+1)=L'(\chi_{-3},-1)
\]
established by C.~Smyth \cite{Sm81} and a circulated conjecture of T.~Chinburg~\cite{Ch84}.
The latter suggests that, given an odd Dirichlet character $\chi_{-N}$, there should be a polynomial with integer coefficients $P_N(x,y)$ for which $L'(\chi_{-N},-1)/\m(P_N)$ is a rational number.%
\footnote{Chinburg also states in \cite{Ch84} a general conjecture for $d$-variable Mahler measures but he shows its truth for $d=1$ only and discusses some evidence for $d=2$.
Unfortunately for the latter, there is an irreparable mistake in his proof of Theorem~2 (namely, the statement of \cite[Lemma~1]{Ch84} is incorrect).}
It seems more plausible to expect the existence of polynomial $P_N(x,y)\in\mathbb Z[x,y]$ such that
\[
\m(P_N)=rL'(\chi_{-N},-1)+\log|s|
\]
for some rational $r\ne0$ and algebraic $s\ne0$.
A statement of this type is shown to be true for the $L$-value $L'(E,0)$ of an elliptic curve $E$ over $\mathbb Q$ with complex multiplication by R.~Pengo \cite[Theorem 4.7]{Pe20}.

\section{Final touch}
\label{s5}

A number $C$ is called an Ap\'ery limit (of order~$d$) if there is an irreducible homogeneous linear difference equation
\[
a_d(n)r(n+d)+a_{d-1}(n)r(n+d-1)+\dots+a_1(n)r(n+1)+a_0(n)r(n)=0
\]
with coefficients $a_d(n),\dots,a_1(n),a_0(n)\in\mathbb Z[n]$, $a_d(n)a_0(n)\ne0$, and its two \emph{rational-valued} solutions $\{q(n)\}_{n\ge n_0}$ and $\{p(n)\}_{n\ge n_0}$ such that $p(n)/q(n)\to C$ as $n\to\infty$. Though we only discuss above the case $d=2$, some constructions are known \cite{CS20,Zu02b} which demonstrate that many `fundamental'  constants, like the values of Riemann's zeta function at positive integers, are Ap\'ery limits.
In analogy with Chinburg's general conjecture from \cite{Ch84}, we would expect that every $L$-value $L(\chi,k)$, for $k\ge2$, associated with a quadratic character $\chi$ is an Ap\'ery limit.

\medskip
\noindent
\textbf{Acknowledgements.}
It is my pleasure to thank Fran\c cois Brunault, Ian Kiming, Nikolay Moshchevitin, Riccardo Pengo, Berend Ringeling and Armin Straub for inspiring conversations and assistance.

A draft of this manuscript (namely, the part in Section~\ref{s3}) was written during my visit in the Mathematics Department of the University of Copenhagen in January 2004.
It is never late to acknowledge the wonderful working conditions I experienced during this visit.
Asmus Schmidt was my host there; our discussions extended my knowledge mathematically and culturally.
One of his favourite (non-mathematical) observations that Asmus constantly stressed on (with proofs!) was a quality of anything marked by the name `Jensen'.
My interest in the Mahler measure grew up later; as Jensen's formula is involved in any single fact about the measure, this serves as my personal motivation (and another evidence for the observation).

\end{document}